\documentclass[11pt]{article}
\usepackage{graphicx}
\usepackage{amsmath}
\usepackage{amsfonts}
\usepackage{epsfig}
\usepackage{setspace}

\newcommand{\be}{\begin{equation}}
\newcommand{\ee}{\end{equation}}
\newcommand{\beqn}{\begin{eqnarray}}
\newcommand{\eeqn}{\end{eqnarray}}
\newcommand{\beqns}{\begin{eqnarray*}}
\newcommand{\eeqns}{\end{eqnarray*}}

\newcommand{\Var}{\mbox{Var}}

\newcommand{\EE}{\ensuremath{{\mathbb E}}}
\newcommand{\II}{\ensuremath{{\mathbb I}}}

\newcommand{\fr}[1]{(\ref{#1})}

\newtheorem{example}{Example}
\newtheorem{lemma}{Lemma}
\newtheorem{theorem}{Theorem}
\newtheorem{corollary}{Corollary}
\newtheorem{remark}{Remark}

\begin{document}

\title{\Large{\bf Nonparametric empirical Bayes estimation based on generalized Laguerre series}}

\author{\large{  Rida Benhaddou}  \footnote{E-mail address: Benhaddo@ohio.edu} \ \ and Matthew A. Connell 
  \\ 
Department of Mathematics, Ohio University, Athens, OH 45701} \date{}

\doublespacing
\maketitle
\begin{abstract}
In this work, we delve into the nonparametric empirical Bayes theory and approximate the classical Bayes estimator by a truncation of the generalized Laguerre series and then estimate its coefficients by minimizing the prior risk of the estimator. The minimization process yields a system of linear equations the size of which is equal to the truncation level. We focus on the empirical Bayes estimation problem when the mixing distribution, and therefore the prior distribution, has a support on the positive real half-line or a subinterval of it. By investigating several common mixing distributions, we develop a strategy on how to select the parameter of the generalized Laguerre function basis so that our estimator {possesses a finite} variance. We show that our generalized Laguerre empirical Bayes approach is asymptotically optimal in the minimax sense. Finally, our convergence rate is compared and contrasted with {several} results from the literature.  \\


{\bf AMS (2000) Subject Classification: 62G05, 62G20, 62G08 }
 \end{abstract} 

\section{Introduction}
 Empirical Bayes (EB) methods are estimation techniques in which the prior distribution, in  the classical Bayesian sense, is estimated from the data. They are powerful tools, in particular when data are generated by repeated execution of the same type of experiment and the size of the dataset is quite large.  Empirical Bayes stands in contrast to classical Bayesian methods, for which the prior distribution is supposed to be fixed before any data are observed. In a typical empirical Bayes setup, one observes two-dimensional random vectors $(x_1, \theta_1), \cdots, (x_N, \theta_N)$, where each $\theta_i$ is distributed according to some unknown prior $g$ and given $\theta_i=\theta$, $x_i$ has the known mixing (conditional) distribution $q(X\mid \theta)$. After the $(N+1)^{st}$ observation $y=x_{N+1}$ is taken, the goal is to estimate its corresponding $\theta_{N+1}$, denoted by $t=\theta_{N+1}$. If the prior {\bf was} known,  Bayes estimator for $t=\theta_{N+1}$ which achieves the smallest mean squared risk would be given by posterior mean 
 \be \label{posm}
 t(y)=\EE\left( \theta\mid x_1, x_2, \cdots, x_N;y\right)=\frac{\int^{\infty}_0\theta q(y\mid \theta)g(\theta)d\theta}{\int^{\infty}_0 q(y\mid \theta)g(\theta)d\theta}=\frac{\Psi(y)}{p(y)}.
  \ee
  Since the prior density $g(\theta)$ is {assumed to be} unknown, we intend to develop an empirical Bayes estimator for \fr{posm}. 
  
  Empirical Bayes estimation methods can be partitioned into two categories: parametric and nonparametric. Parametric empirical Bayes procedures rely on the knowledge of the parametric form of the prior distribution. In particular, if the mixing distribution and its prior are assumed to take on simple parametric forms, such as the case of mixing distributions with simple conjugate priors, then the empirical Bayes problem becomes that of estimating the marginal $p(y)$ using a set of empirical quantities. One way is to approximate the marginal $p(y)$ by the means of the {maximum likelihood method}, method of moments or expectation-maximization (EM) algorithm, using past data, which makes it easier to express the parameters of the prior in terms of their empirical counterparts. The second step would be to replace the empirical hyper-parameters of the prior in the right-hand side of \fr{posm}. On the other hand, in nonparametric empirical Bayes methods, the prior distribution and its parametric form are not specified at all. One of the techniques to nonparametric empirical Bayes methods is to estimate the top {$\Psi(y)$} and bottom {$p(y)$} of estimator \fr{posm} separately and {\bf then} compute the ratio{, this includes Singh~(1976, 1979), Nogami~(1988), Tiwari and Zalkikar~(1990), Ma and Balalrishnan~(2000) among others}. The other set of methods rely on estimating the ratio in estimator \fr{posm} directly {and this one includes Robbins~(1983), Pensky and Ni~(2000), Pensky and Alotaibi~(2005) and Benhaddou and Pensky~(2013), among others}. 
          
 Since the seminal work of Robbins~(1955,\ 1964) there has been a great array of papers that address EB estimation, and the list includes Efron and Morris~(1977), Morris~(1983), Louis~(1984) and Casella~(1985) in the parametric EB case, and Singh~(1976, 1979), Nogami~(1988), Tiwari and Zalkikar~(1990), Datta~(1991, 2000), Walter and Hamedani~(1991), Pensky~(1997a, 1997b, 2002), Ma and Balakrishnan~(2000), Brown and Greenshtein~(2009) and Benhaddou and Pensky~(2013)  in the nonparametric EB case. In 1983, Herbert Robbins introduced the linear empirical Bayes estimation procedure which relies on approximating the Bayes rule $t(y)$ locally by a linear function of $y$ and estimating its coefficients by minimizing the prior risk of the estimator. This technique is very efficient computationally and it was immediately put into practical use  (see e.g., Ghosh and Meeden~(1986), Ghosh and Lahiri~(1987)). The shortcoming of a linear EB estimator is that it possesses a large bias and it is optimal only in the class of estimators that are linear in $y$. Pensky and Ni~(2000) proposed an extension of the linear EB estimator of Robbins~(1983) to approximating $t(y)$ by an algebraic polynomial, and Pensky and Alotaibi ~(2005) and Benhaddou and Pensky~(2013) devised an estimator based on approximating the Bayes rule $t(y)$ via wavelet series. In the latter wavelet-based approaches, the expansion of  $t(y)$ is carried out over scaling functions at resolution level $m$ and the coefficients of the expansion are estimated by minimizing the prior risk, and in Benhaddou and Pensky~(2013) the choice of resolution level $m$ is performed via Lepski~(1997) method.  {In this paper, we propose a series expansion approach based on a system of generalized Laguerre functions basis}. 

    The application of Laguerre series to nonparametric estimation has become more popular as of late and the list includes the application to density estimation in Comte et al.~(2008), Comte and Genon-Catalot~(2015) and Dussap~(2021), the estimation of linear functionals of a density function in  Mabon~(2016), Laplace deconvolution in Vareschi~(2015), Comte, Cuenod, Pensky, and Rozenholc~(2017) and Benhaddou, Pensky and Rajapakshage~(2019), and to nonparametric regression estimation in Benhaddou~(2021). 
    
In the present work, we extend the linear empirical Bayes estimator of Robbins~(1983) by approximating the right-hand side of \fr{posm} by a truncation of the generalized Laguerre series and estimating its coefficients by minimizing the prior risk.  The minimization leads to a system of linear equations that is well-conditioned thanks to the interesting features of the generalized Laguerre function basis. We only consider the EB estimation problem when the mixing distribution $q(X\mid \theta)$ is defined on the positive real half-line or a subinterval of it. It is demonstrated that the choice of the truncation level depends on the balance between the variance and the bias of the estimator and if chosen properly the proposed EB estimator attains asymptotically optimal convergence rates in the minimax point of view.   In addition, by investigating several common mixing distributions, we develop a strategy on how to select the parameter of the generalized Laguerre function basis so that our estimator  {possesses a finite} variance. Our results are then compared to several similar existing nonparametric procedures found in the literature. 

\subsection{Layout of the paper}

The rest of the paper is organized as follows. Section $2$ describes in detail the Laguerre-based empirical Bayes estimation algorithm, explores several examples of common mixing distributions and motivates the particular choice of the  parameter of  generalized Laguerre function system to be used in each case. Section $3$ presents the main results of the paper, that is, it provides asymptotic results on the variance and bias components, discusses the choice of the truncation level of the Laguerre-based empirical Bayes estimator to be used, and then gives the asymptotic assessment of the estimation error and the minimax optimality of the truncated Laguerre series procedure. Finally, an {appendix} contains supplementary results and the proofs of  our theoretical statements that appear throughout the paper. 

\subsection{Notation}

 For the rest of the paper, let $\|h\|$ denote the vector norm of the vector $h$. Given a matrix ${\bf A}$, let ${\bf A}^T$ be its transpose, $\lambda_{\max}({\bf A})$ be its largest eigenvalue in magnitude, $\|{\bf A}\|_F=\sqrt{Tr\left(A^TA\right)}$ and $\|{\bf A}\|_{sp}=\lambda_{\max}\left({\bf A}^T{\bf A}\right)$ be, respectively, its Frobenius and the spectral norms. In addition, let $(a\vee b)=\max(a, b)$ and $(a\wedge b)=\min(a, b)$. Finally, {for the sequences $\{a_n\}$ and $\{b_n\}$ of positive real numbers, let $a_n\asymp b_n$ denote the property that there exist positive constants $c_1$ and $c_2$ independent of $n$ such that $c_1 \leq a_n/b_n\leq c_2$, with $c_1< c_2< \infty$}.

 \section{Estimation Algorithm}
 
 In order to construct an empirical Bayes estimator for $t(y)$ defined above,  consider the orthonormal basis that consists of the system of generalized Laguerre functions 
\be
\varphi^{(a)}_k(x)={\left[\frac{k!}{\Gamma(k+a+1)}\right]^{1/2}} e^{-x/2}x^{a/2}L^{(a)}_k(x), \ \ k=0, 1, \cdots, 
\ee
where $L^{(a)}_k(t)$ are generalized Laguerre polynomials with parameter $a$, $a\geq 0$ (see. e.g., Gradshtein and Ryzhik~(1980), Section 8.97). Approximate $t(y)$ by the truncated Laguerre series, for some {relatively large} $M$, 
\be \label{lagmattr}
t_M(y)= \sum^{M-1}_{l=0}\theta_l\varphi^{(a)}_l(y),
\ee
and estimate the coefficients $\theta_l$, $l=0, 1, 2, \cdots, M-1$, by minimizing the integrated mean squared error  
 \be \label{gjk}
 \int^{\infty}_0\int^{\infty}_0\left[t_M(y)-z\right]^2q(y\mid z)g(z)dzdy.
\ee 
The first order necessary condition for minimizing \fr{gjk} with respect to $\theta_l$, $l=0, 1, 2, \cdots, M-1$, yields the system of linear equations 
\be\label{est-f}
{\bf A}_M{\pmb{ \Theta}_{M}}={\bf C}_M,
\ee \label{thetest}
where ${\bf A}_M$ is the $M\times M$ matrix with elements 
\be \label{alk}
A_{lk}=\int^{\infty}_{0}\varphi^{(a)}_l(x)\varphi^{(a)}_k(x)p(x)dx, \ \ l, k=0, 1, 2, \cdots, M-1,
\ee
 and ${\bf C}_M$ is the $M$-dimensional vector with elements 
\be \label{ck}
c_{k}=\int^{\infty}_{0} \varphi^{(a)}_k(x)\Psi(x)dx,\ \  k=0, 1, 2, \cdots, M-1.
\ee
The entries of the matrix ${\bf A}_M$ are unknown but can be estimated from the data via the sample {averages}  
\be \label{ahat}
\widehat{A}_{lk}=N^{-1}\sum^N_{i=1}\varphi^{(a)}_l(X_i)\varphi^{(a)}_k(X_i), \ \ {l, k=0, 1, 2, \cdots, M-1}.
\ee
{As for the entries $c_{k}$ of the vector ${\bf C}_M$, notice that the right-hand side of equation \fr{ck} is not expressed in the form of marginal expectation as is in the case of equation \fr{alk}, that is, it does not involve the marginal distribution $p(x)$ and therefore it cannot be estimated directly from the data, $X_1, X_2, \cdots, X_N$. As a result, one needs to find equivalent expression for $c_{k}$ that involves the marginal distribution $p(x)$}. {This can be achieved by finding} the functions $U_k(x)$, $k=0, 1, 2, \cdots, M-1$, such that for any $\theta>0$ 
\be \label{uk}
\int^{\infty}_0q(x\mid \theta)U_k(x)dx=\int^{\infty}_0\theta q(x\mid \theta)\varphi^{(a)}_k(x)dx.
\ee
{Observe} the connection between $U_k(x)$ in equation \fr{uk} and the entries $c_{k}$ in \fr{ck}, the expectation of $U_k(x)$ over the marginal distribution $p(x)$ is 
\beqn \label{U-def}
\EE[U_k(x)]=\int^{\infty}_0U_k(x)p(x)dx=\int^{\infty}_0\left[\int^{\infty}_0U_k(x)q(x\mid \theta)dx\right]g(\theta)d\theta=\int^{\infty}_0\varphi^{(a)}_k(x)\Psi(x)dx.
\eeqn
{It is worth pointing out that} the construction of the functions $U_k(x)$ is possible for a number of known mixing distributions that are defined on the positive real-line or at least a subinterval of it. Once such functions are derived, then entries in \fr{ck} can be estimated by {their sample averages}
\be \label{chat}
\widehat{c}_k=N^{-1}\sum^N_{i=1}U_k(X_i), \ \  {k=0, 1, 2, \cdots, M-1}.
\ee
 \subsection{Estimation of generalized Laguerre coefficients}
 
Once, entries of matrix ${\bf A}_M$ and vector ${{\bf C}}_M$ are estimated, system \fr{est-f} may then be replaced by the surrogate 
\be\label{est-h}
{\bf \widehat{A}}_M \widehat{{\pmb{{ \Theta}}}}_M=\widehat{{{\bf C}}}_M.
\ee 
However, since ${\bf \widehat{A}}_M$ and $\widehat{{{\bf C}}}_M$ are asymptotically normal, estimator $\widehat{{\pmb{{\Theta}}}}_M=\left({\bf \widehat{A}}_M\right)^{-1}\widehat{{{\bf C}}}_M$ may not have a finite expectation. To guarantee that the estimator of ${{\pmb{{\Theta}}}}_M$ has finite expectation, we choose $\delta(N) >0$ and consider instead an estimator of the form 
\be\label{est-theta}
 \widehat{{\pmb{{\Theta}}}^{\delta}}_M=\left({\bf \widehat{A}}_M+\delta(N){\bf I}_M\right)^{-1} \widehat{{{\bf C}}}_M,
\ee 
where ${\bf I}_M$ is the identity matrix of size $M$ and $\delta(N)$ is a positive quantity. Notice that now the matrix $\left({\bf \widehat{A}}_M+\delta(N){\bf I}_M\right)$ is positive definite, and therefore, it is invertible. Thus, consider the generalized Laguerre empirical Bayes estimator  
\be \label{lagmat-hat}
\widehat{t}_M(y)= \sum^{M-1}_{l=0} (\widehat{{\pmb{{\Theta}}}^{\delta}}_M)_l\varphi^{(a)}_l(y),
\ee
{where the choices of $M$ and $\delta(N)$ depend on the sample size $N$ and will be determined later}. \\
Next, we will take a look at a variety of common mixing probability density functions to see how such functions $U_k(x)$ can be found. 

 \subsection{Examples and derivation of the functions $U_k(x)$}
 To see how the generalized Laguerre procedure works, how to find the functions $U_k(x)$ and how to decide on the parameter $a$ of generalized Laguerre function basis to use, we will investigate several known mixing distributions defined on $(0, \infty)$ or some subinterval of it. The desired feature to look for in the functions $U_k(x)$ is to have finite $L^p$-norm, with $p \geq 1$, and the parameter $a$ of generalized Laguerre function basis may be chosen specifically to guarantee that. 
   
\begin{example} \label{ex1}
Let $q(x\mid \theta)$ be a uniform distribution. In particular, let 
\be \label{unif-con}
q(x\mid \theta)=\frac{1}{\theta} \II\left(0< x< \theta\right),  a\leq \theta \leq b.
\ee
It is worth pointing out that uniform distribution of the form \fr{unif-con} was considered in Nogami~(1988) where a kernel-based empirical Bayes procedure was suggested to estimate $\theta$. 
Notice that by \fr{uk}, we have
\be \label{uniu}
\int^{\theta}_0U_l(x)dx=\theta \int^{\theta}_0\varphi^{(a)}_l(x)dx.
\ee
Then, differentiating both sides of \fr{uniu} with respect to $\theta$ yields
\beqn \label{uuxni}
U_l(\theta)=\int^{\theta}_0\varphi^{(a)}_l(x)dx+\theta \varphi^{(a)}_l(\theta).
\eeqn
Finally, replacing $\theta$ with $x$ in \fr{uuxni}, gives
\be \label{uuni}
U_l(x)=\int^{x}_0\varphi^{(a)}_l(z)dz+x \varphi^{(a)}_l(x).
\ee
Now, take the parameter of the generalized Laguerre function basis $a=0$ and let us find an upper bound for the quadratic norm of $U_l(x)$. Since $a\leq \theta \leq b$, then $a\leq x \leq b$ and therefore, for any $l\geq0$, one has 
\beqn
\int^{b}_aU^2_l(x)dx \leq 2 \int^b_ax^2 \varphi^2_l(x)dx+2 \int^b_a\left(\int^x_0\varphi_l(z)dz\right)^2dx\leq c.
\eeqn
\end{example}
\begin{example} \label{ex2}
Let $q(x\mid \theta)$ be a Pareto distribution with scale parameter $\theta >0$, and shape parameter $\alpha >2$. In particular, let 
\be \label{paret-con}
q(x\mid \theta)=\frac{\alpha \theta^{\alpha}}{x^{\alpha +1}},  \ x \geq \theta\ and\ 0 < \theta\leq \theta_o.
\ee
 Pareto mixing of form \fr{paret-con} was considered in Tiwari and Zalkikar~(1990) where, in the spirit of Nogami~(1988), a kernel-based empirical Bayes approach was proposed  to estimate $\theta$.  Observe that by \fr{uk}, we have
\be \label{eq-par-p}
\int^{\infty}_{\theta}x^{-\alpha -1}U_l(x)dx=\theta \int^{\infty}_{\theta}x^{-\alpha-1}\varphi^{(a)}_l(x)dx.
\ee
Now, differentiating both sides of \fr{eq-par-p} with respect to $\theta$ yields
\beqn
-\frac{1}{\theta^{\alpha+1}}U_l(\theta)=\int^{\infty}_{\theta}\frac{1}{x^{\alpha+1}}\varphi^{(a)}_l(x)dx-\theta^{-\alpha} \varphi^{(a)}_l(\theta).
\eeqn
Finally, rearranging and replacing $\theta$ with $x$ gives
\be \label{uuni}
U_l(x)=-x^{\alpha+1}\int^{\infty}_x\frac{1}{z^{\alpha+1}}\varphi^{(a)}_l(z)dz+x \varphi^{(a)}_l(x).
\ee
Here, take the parameter of the generalized Laguerre function basis $a=0$ and let us find an upper bound for the quadratic norm of $U_l(x)$. Since $x \geq \theta$ and $0< \theta \leq \theta_o$, then $\theta_o\leq x$ and therefore, using equation $2.5$ in Muckenhoupt~(1970) and for $\theta_o>1$ and any $l\geq0$, yields 
\beqn
\int^{\infty}_{\theta_o}U^2_l(x)dx \leq 2 \int^{\infty}_{\theta_o}x^2 \varphi^2_l(x)dx+2 \int^{\infty}_{\theta_o}x^{2\alpha+2}\left(\int^{\infty}_x\frac{1}{z^{\alpha+1}}\varphi_l(z)dz\right)^2dx\leq c(l^{7/3}\vee 1).
\eeqn
\end{example}
\begin{example} \label{ex3}
Let $q(x\mid \theta)$ be a beta distribution with unknown shape parameter $\theta >0$, and known shape parameter $\alpha >0$. In particular, let 
\be \label{var-bias-b}
q(x\mid \theta)=\frac{\Gamma(\theta + \alpha)}{\Gamma(\theta)\Gamma(\alpha)}x^{\alpha-1}(1-x)^{\theta-1},  \ x \in (0, 1).
\ee
Notice that by \fr{uk}, one has 
\be \label{eq-par-b}
\int^{1}_{0}x^{\alpha -1}(1-x)^{\theta-1}U_l(x)dx=\theta \int^{1}_{0}x^{\alpha-1}(1-x)^{\theta-1}\varphi^{(a)}_l(x)dx.
\ee
Now, integrating by parts the right-hand side of \fr{eq-par-b}  yields
\beqn
\int^{1}_{0}x^{\alpha -1}(1-x)^{\theta-1}U_l(x)dx=\int^{1}_{0}(1-x)^{\theta}\left[(\alpha-1)x^{\alpha-2}\varphi^{(a)}_l(x)+x^{\alpha-1}\left(\varphi^{(a)}_l(x)\right)'\right]dx.
\eeqn
Finally, rearranging gives
\be \label{ubet}
U_l(x)=(\alpha-1)(1-x)x^{-1}\varphi^{(a)}_l(x)+(1-x) \left(\varphi^{(a)}_l(x)\right)'.
\ee
Note that to ensure that \fr{ubet} has a finite $L^2$-norm, take the parameter of the generalized Laguerre function basis $a=2$. Let us now find an upper bound for the norm of $U_l(x)$. Indeed, by equation $2.5$ in Muckenhoupt~(1970) and for any $l\geq0$, computing the norm yields 
\beqn
\int^{1}_{0}U^2_l(x)dx \leq 2(\alpha-1)^2 \int^{1}_{0} x^{-2}\left(\varphi^{(2)}_l(x)\right)^2dx+2 \int^{1}_{0}\left(\left[\varphi^{(2)}_l(x)\right]'\right)^2dx\leq c(l\vee 1).
\eeqn
\end{example}
\begin{example} \label{ex4}
Let $q(x\mid \theta)$ be an exponential distribution with  scale parameter $\theta >0$. In particular, let 
\be \label{var-bias-b}
q(x\mid \theta)=\theta \exp\{-x\theta\},  \ x \in [0, \infty).
\ee
Observe that by \fr{uk}, we  have
\be \label{eq-par-ex}
\int^{\infty}_{0} \exp\{-x\theta\}U_l(x)dx=\theta \int^{\infty}_{0} \exp\{-x\theta\}\varphi^{(a)}_l(x)dx.
\ee
Now, for any $a >0$,  integrating by parts the right-hand side of \fr{eq-par-ex}  yields
\beqn
\int^{\infty}_{0} \exp\{-x\theta\}U_l(x)dx=\int^{\infty}_{0}\exp\{-x\theta\}\left[\left(\varphi^{(a)}_l(x)\right)'\right]dx.
\eeqn
Finally, rearranging and solving gives
\be \label{u-ex}
U_l(x)= \left(\varphi^{(a)}_l(x)\right)'.
\ee
Remark that to ensure that \fr{u-ex} has a finite $L^2$-norm, we take the parameter of the generalized Laguerre function basis $a=2$. Thus, by equation $2.2$ in Muckenhoupt~(1970) and for any $l\geq0$, computing the norm of $U_l(x)$ yields  
\beqn
\int^{\infty}_{0}U^2_l(x)dx = \int^{\infty}_{0}\left(\left[\varphi^{(2)}_l(x)\right]'\right)^2dx\leq c(l\vee 1).
\eeqn
\end{example}
\begin{example} \label{ex5}
Let $q(x\mid \theta)$ be a Rayleigh distribution with scale parameter $\theta >0$. In particular, let 
\be \label{var-bias-r}
q(x\mid \theta)=\theta x \exp\{-x^2\theta/2\},  \ x \in [0, \infty).
\ee
Notice that by \fr{uk}, we have
\be \label{eq-par-ray}
\int^{\infty}_{0} x\exp\{-x^2\theta/2\}U_l(x)dx=\theta \int^{\infty}_{0}x \exp\{-x^2\theta/2\}\varphi^{(a)}_l(x)dx.
\ee
Now, for any $a >0$, integrating by parts the right-hand side of \fr{eq-par-ray}  yields
\beqn
\int^{\infty}_{0} x\exp\{-x^2\theta/2\}U_l(x)dx=\int^{\infty}_{0}\exp\{-x^2\theta/2\}\left[\left(\varphi^{(a)}_l(x)\right)'\right]dx.
\eeqn
Finally, rearranging and solving gives
\be \label{u-ray}
U_l(x)=\frac{1}{x} \left(\varphi^{(a)}_l(x)\right)'.
\ee
In this case, to ensure that \fr{u-ray} has a finite $L^2$-norm, we take the parameter of the generalized Laguerre function basis $a=4$. Therefore, by equation $2.2$ in Muckenhoupt~(1970) and for any $l\geq0$, computing the norm of $U_l(x)$ yields
\beqn
\int^{\infty}_{0}U^2_l(x)dx = \int^{\infty}_{0}\left(\frac{1}{x}\left[\varphi^{(4)}_l(x)\right]'\right)^2dx\leq c(l\vee 1).
\eeqn
\end{example}
\begin{example} \label{ex6}
Let $q(x\mid \theta)$ be a Weibull distribution with unknown scale parameter $\theta >0$ and known shape parameter $\alpha >0$. In particular, let 
\be \label{var-bias-w}
q(x\mid \theta)=\alpha \theta x^{\alpha-1} \exp\{-x^{\alpha}\theta\},  \ x \in [0, \infty).
\ee
Notice that by \fr{uk}, we obtain 
\be \label{eq-par-w}
\int^{\infty}_{0} x^{\alpha-1}\exp\{-x^{\alpha}\theta\}U_l(x)dx=\theta \int^{\infty}_{0}x^{\alpha-1} \exp\{-x^{\alpha}\theta\}\varphi^{(a)}_l(x)dx.
\ee
Now, for any $a >0$, integrating by parts the right-hand side of \fr{eq-par-w} yields
\beqn
\int^{\infty}_{0} x^{\alpha-1}\exp\{-x^{\alpha}\theta\}U_l(x)dx=\frac{1}{\alpha}\int^{\infty}_{0}\exp\{-x^{\alpha}\theta\}\left[\left(\varphi^{(a)}_l(x)\right)'\right]dx.
\eeqn
Finally, rearranging and solving gives
\be \label{u-w}
U_l(x)=\frac{1}{\alpha x^{\alpha-1}} \left(\varphi^{(a)}_l(x)\right)'.
\ee
In this case, to ensure that \fr{u-w} has a finite $L^2$-norm, take the parameter of the generalized Laguerre function basis $a=2\alpha$. Hence, by equation $2.2$ in Muckenhoupt~(1970) and for $l\geq0$, computing the $L^2$-norm of $U_l(x)$ yields 
\beqn
\int^{\infty}_{0}U^2_l(x)dx = \int^{\infty}_{0}\left(\frac{1}{\alpha x^{\alpha-1}}\left[\varphi^{(2\alpha)}_l(x)\right]'\right)^2dx\leq c(l\vee 1).
\eeqn
\end{example}
 \begin{remark}	
 Note that in all of the examples above, the choice of parameter of the generalized Laguerre function basis, $a$, guarantees that the functions $U_l(x)$ will have finite $L^p$-norm for any  $p \geq 1$, not just the $L^2$-norm. 
  \end{remark} 
   \section{Minimax convergence rates and adaptivity}
  In order to evaluate the accuracy of the proposed Laguerre-based empirical Bayes  estimator, we need to decide what risk measure to use.  In general, the goodness of an empirical Bayes estimator $\widehat{t}_N(y)$ can be assessed by its posterior risk 
   \be \label{postrisk}
      R(\widehat{t}_N(y))=\frac{1}{p(y)}\EE\int^{\infty}_0\left(\widehat{t}_N(y)-\theta\right)^2q(y\mid \theta)g(\theta)d\theta .  
   \ee
   The use of the posterior risk of form \fr{postrisk} enables us to assess the minimax optimality of an empirical Bayes estimator of the majority of mixing distributions via the comparison of the upper bounds of the empirical Bayes estimator to the lower bounds for the risk derived in Benhaddou and Pensky~(2013). Observe that \fr{postrisk} can be partitioned as follows 
      \beqn
   R(\widehat{t}_N(y)) &=& \EE\left[\widehat{t}_N(y)-t(y)\right]^2 + \frac{1}{p(y)}\int^{\infty}_0\left({t}(y)-\theta\right)^2q(y\mid \theta)g(\theta)d\theta. \label{postrisk1}
               \eeqn
  Notice now that the second term in the right-hand side of \fr{postrisk1} represents the pure posterior risk of the classical Bayes estimator of equation \fr{posm}, and therefore it is irrelevant to the assessment of the precision of our Laguerre-based empirical Bayes procedure. So, it makes sense to measure the performance of our estimator \fr{lagmat-hat} by the first term only, denote it by $R_N(y)$. As a result, we will evaluate the proposed estimator based on the risk
     \be \label{rn}
   R_N(y)=\EE\left[\widehat{t}_M(y)-t(y)\right]^2.
   \ee
 Now, observe that \fr{rn} itself can be partitioned into the usual sum of the random error component, or variance term, and the systematic error component, or bias term, as follows
 \be \label{rn12}
   R_N(y)\leq 2R_1(y)+2R_2(y),
   \ee
where
\beqn
R_1(y)&=&\left[t_M(y)-t(y)\right]^2, \label{r1}\\
R_2(y)&=& \EE\left[\widehat{t}_M(y)-t_M(y)\right]^2 \label{r2}.
\eeqn
Since we are using the posterior risk as the measure of the performance of the empirical Bayes estimator \fr{lagmat-hat}, we will treat $y$ as a fixed quantity.
 \subsection{Asymptotic evaluation of the bias term}
In order to evaluate the systematic error term \fr{r1}, define the set $\pmb{\Omega}_y(M)$ as follows 
\be \label{omset}
\pmb{\Omega}_{y}(M)=\left\{x: |x-y|\leq \frac{s}{M}\right\}.
\ee
Then, applying the substitution ${z}{M^{-1}}=x-y$ to the right-hand sides of \fr{alk}, and then Maclaurin series expansion to $p(y+zM^{-1})$ in \fr{alk} yields
\beqn
A_{lk}&=&p(y)\int^{\infty}_0\varphi_l(u)\varphi_k(u)du +M^{-1}p'(y)\int^{\infty}_0\varphi_l(u)\varphi_k(u)(u-y)du+ \frac{1}{2}M^{-2}p''(y)\int^{\infty}_0\varphi_l(u)\varphi_k(u)(u-y)^2du\nonumber\\
&+&\cdots + \frac{1}{r!}M^{-r}p^{(r)}(y)\int^{\infty}_0\varphi_l(u)\varphi_k(u)(u-y)^rdu+o(M^{-r}).
\eeqn
Hence, the matrix ${\bf A}_M$ takes the form 
\beqn \label{amser}
{\bf A}_M=p(y)I_M+\sum^r_{h=1}M^{-h}\frac{p^{(h)}(y)}{h!} {\bf \Phi}^{(h)}_M+o(M^{-r}),
\eeqn
where ${\bf \Phi}^{(h)}_M$ are $M\times M$ matrices with elements 
\beqns
\int^{\infty}_{0}\varphi_l(u)\varphi_k(u)(u-y)^hdu.
\eeqns
In addition, applying the substitution ${z}{M^{-1}}=x-y$ to the right-hand sides of \fr{ck}, and then Maclaurin series expansion to $\Psi(y+zM^{-1})$ in\fr{ck} yields
\beqn
c_{k}&=&\Psi(y)\int^{\infty}_0\varphi_k(u)du +M^{-1}\Psi'(y)\int^{\infty}_0\varphi_k(u)(u-y)du+ \frac{1}{2}M^{-2}\Psi''(y)\int^{\infty}_0\varphi_k(u)(u-y)^2du\nonumber\\
&+&\cdots + \frac{1}{r!}M^{-r}\Psi^{(r)}(y)\int^{\infty}_0\varphi_k(u)(u-y)^rdu+o(M^{-r}).
\eeqn
Therefore, the vector ${\bf C}_M$ takes the form
\beqn \label{cmser}
{\bf C}_M=\sum^r_{h=0}M^{-h}\frac{\Psi^{(h)}(y)}{h!}\Lambda^{(h)}_M+o(M^{-r}),
\eeqn
where $\Lambda^{(h)}_M$ are $M$-dimensional vectors with elements 
\beqns
\int^{\infty}_{0}\varphi_k(u)(u-y)^hdu.
\eeqns
Notice that as $M\rightarrow \infty$, $A_{lk} \rightarrow p(y)\int^{\infty}_0\varphi_l(u)\varphi_k(u)du$ and $c_k \rightarrow \Psi(y)\int^{\infty}_0\varphi_k(u)du$. Hence, by Theorem $12$ of Muckenhoupt~(1970), which guarantees the mean convergence of the partial sum of the Laguerre polynomial series, as $M\rightarrow \infty$
\be \label{lag-asymp}
t_M(y)= \sum^{M-1}_{l=0}({\bf A}_M^{-1}{\bf C}_M)_l\varphi^{(a)}_l(y) \rightarrow \frac{\Psi(y)}{p(y)}=t(y).
\ee
Therefore, the system of equations \fr{est-f} is well-conditioned and {\bf the} following statement is true.
  \begin{lemma} \label{lem:Bias}
Let $p(x)$ and $\Psi(x)$ be $r$ times continuously differentiable in the neighborhood $\pmb{\Omega}_y$ of $y$ and let $\pmb{\Omega}_{y}(M)$ be a subset of  $\pmb{\Omega}_y$. Let $r^*=(r\wedge s)$, where $s$ is a positive real number such that $s >1$. Then, for any $y$ such that $p(y)\neq 0$, as $M \rightarrow \infty$, one has 
\be \label{var-bias}
R_1(y)=o\left(M^{-2r^*}\right).
\ee
\end{lemma}
 \subsection{Asymptotic evaluation of the variance term}
 We assume that the functions $U_l(x)$ defined in \fr{U-def} satisfy the following condition. \\
  \noindent
{\bf Assumption A.2.}  There exist constants $c_u>0$ and $\beta \geq 0$, independent of $N$ and $l$, such that
\be
\int^{\infty}_{0}U^2_l(x)dx\leq c_u (l^{\beta}\vee 1). \label{LRD}
\ee
 \begin{remark}	
 {\bf Assumption A.2}  is {\bf inspired by the behavior of the functions $U_l(x)$} under mixing distributions discussed in {\bf Examples \ref{ex1}} through {\bf \ref{ex6}}. 
  \end{remark} 
  \begin{lemma} \label{lem:var}
Let ${\bf A}_M$, ${{\bf C}}_M$,  ${\bf \widehat{A}}_M$ and $\widehat{{{\bf C}}}_M$ be defined in \fr{alk}, \fr{ck}, \fr{ahat} and \fr{chat}, respectively, and  let condition \fr{LRD} hold. Then,  as $N, M \rightarrow \infty$, one has 
\beqn
\EE \| \widehat{\bf A}_M-{\bf A}_M\|_{sp}^2&=&O\left(M^{2}N^{-1}\right). \label{var-a}\\
\EE \| \widehat{{\bf C}}_M-{{\bf C}}_M\|^2&=&O\left(M^{\beta+1}N^{-1}\right). \label{var-c}
\eeqn
In addition,
\beqn
\EE \| \widehat{\bf A}_M-{\bf A}_M\|_{sp}^4&=&O\left(M^{2}N^{-2}\right). \label{var-a2}\\
\EE \| \widehat{{\bf C}}_M-{{\bf C}}_M\|^4&=&O\left(M^{2\beta+1}N^{-2}+ MN^{-3}\max_{l\leq M-1}\int^{\infty}_{0}U^4_l(x)dx\right). \label{var-c2}
\eeqn\end{lemma}
  \begin{lemma} \label{lem:largd}
Let ${\bf A}_M$, ${{\bf C}}_M$,  ${\bf \widehat{A}}_M$ and $\widehat{{{\bf C}}}_M$ be defined in \fr{alk}, \fr{ck}, \fr{ahat} and \fr{chat}, respectively. Then, for any $\gamma > 0$, one has 
\beqn
\Pr\left( \| \widehat{\bf A}_M-{\bf A}_M\|_F^2 > M^2 \gamma^2 N^{-1}\ln(N) \right)& \leq&2 M^2N^{-\tau^2}, \label{devi-a}
\eeqn
where $\tau^2=\frac{\gamma^2}{8\|p\|_{\infty}\|\varphi_k\|_{\infty}}$.
\end{lemma}
 \begin{lemma} \label{lem:R2}
Let $\delta^2(N)=N^{-1}M^2$ and choose $\gamma$ in \fr{devi-a} such that $\gamma^2=\frac{N\|A^{-1}\|_{sp}^2}{4M^2\ln(N)}$. Then, as $M, N \rightarrow \infty$, one has 
\be \label{R-var}
 R_2(y)=O\left({N^{-1}M^{(\beta\vee 1)+1}}\right),
\ee
provided that $N^{-1}M^2 \rightarrow 0$, as $N \rightarrow \infty$ and $M\max_{l\leq M-1}\int^{\infty}_0U^4_l(x)dx=o(N^3)$.
\end{lemma}
 \subsection{Convergence rate and minimax optimality}
{Notice that  the variance component \fr{R-var} and the bias component \fr{var-bias} behave differently to changes in $M$. That is, while the variance is an increasing function of $M$, the bias is a decreasing function of $M$ and therefore choosing an $M$ that is too small or too large will affect the risk in \fr{rn} in a detrimental way}. In order to minimize the posterior risk in \fr{rn}, we choose the truncation level $M$ such that the two errors in \fr{var-bias} and \fr{R-var} are balanced.  According to {\bf Lemmas} {\bf \ref{lem:Bias}} and {\bf \ref{lem:R2}} the posterior risk \fr{rn} is at its minimum when the truncation level $M=M_o$ is such that, as $N \rightarrow \infty$,
\be \label{R-min}
M_o\asymp N^{\frac{1}{(\beta\vee 1)+2r^*+1}}.
\ee
The following theorem gives the asymptotic upper bound for the posterior risk of estimator \fr{lagmat-hat} with truncation level \fr{R-min}. 
\begin{theorem} \label{th:upperbds-2}
Let $p(x)$ and $\Psi(x)$ be $r$ times continuously differentiable in the neighborhood $\pmb{\Omega}_y$ of $y$ and let $\pmb{\Omega}_{y}(M)$  be a subset of $\pmb{\Omega}_y$. Let $\widehat{t}_M(y)$ be the Laguerre estimator defined in \fr{lagmat-hat} with $M$ given in \fr{R-min} and $\delta(N)$ in \fr{est-theta} chosen as $\delta^2(N)=N^{-1}M^2$.  Then, under the conditions of {\bf Lemma \ref{lem:Bias}} and {\bf Lemma \ref{lem:R2}}, as $N \rightarrow \infty$, the posterior risk is of the order 
 \be \label{upperbds-2}
  R_N(y)=O\left( N^{-\frac{2r^*}{(\beta\vee 1)+2r^*+1}}\right).
   \ee
\end{theorem}
 \begin{remark}	
 Lower bounds for the risk of the nonparametric empirical Bayes estimators were {\bf established} in Benhaddou and Pensky~(2013). According to their equation $(4.5)$, with their quantity $r_1$ chosen as $r_1=r^*+1/2(\beta\vee 1)$ and our own {\bf Theorem \ref{th:upperbds-2}}, estimator \fr{lagmat-hat} with $M$ chosen according to \fr{R-min} is  asymptotically optimal in the minimax sense.
  \end{remark} 
  \begin{remark}	  
The convergence rate is expressed in terms of the parameter $\beta$, which is associated with the mixing distribution $q(x\mid \theta)$ under consideration and the choice of the generalized Laguerre basis, as well as the parameter $r$ associated with the smoothness of the functions $p(x)$ and $\Psi(x)$.  Nevertheless, selecting an appropriate parameter $a$  for the generalized Laguerre basis has a deciding effect on the value of $\beta$.
  \end{remark}
 \begin{remark}	  
{Singh~(1979) considered the empirical Bayes estimation when the mixing distribution belongs to the general one-parameter exponential family based on the kernel method, and was able to obtain convergence rates of order $O\left(N^{-\frac{2(r-1)}{2r+1}}\right)$, in their notation. This implies that our rates match theirs when $3r^*+(\beta\vee 1)+1=r((\beta\vee 1)+1)$. However, our rates are faster than theirs whenever $r < \frac{3r^*}{(\beta\vee 1)+1}+1$, and slower otherwise}. 
  \end{remark}
    \begin{remark}	
The case of  a uniform mixing distribution was considered in Nogami~(1988) where a kernel-based empirical Bayes procedure was proposed to estimate $\theta$, and under certain conditions on the prior, it was shown to achieve a convergence rate of order $O\left( N^{-\frac{1}{2}}\right)$. According to {\bf Example \ref{ex1}}, our Laguerre-based estimator with the choice of parameter $a=0$ generalized Laguerre function basis is appropriate. In such case, $\beta=0$ and thus, the minimax risk will be of order $O\left(N^{-\frac{r^*}{r^*+1}}\right)$. Notice that our rate matches that of Nogami~(1988) for the specific case of $r^*=1$, but {it is} faster than that of Nogami~(1988) for $r^*>1$. In addition, for $r^*=r$ and under the same conditions on $p(x)$ and $\Psi(x)$, our rate is slightly slower  than {that} in Benhaddou and Pensky~(2013) in their treatment of problem based on the wavelet approach for this particular mixing distribution (their rate is $N^{-\frac{2r}{2r+1}}$).
  \end{remark}
  \begin{remark}	
The case of  a Pareto mixing distribution was considered in Tiwari and Zalkikar~(1990) where, in the spirit of Nogami~(1988), a kernel-based empirical Bayes approach was used  to estimate $\theta$, and under certain conditions on the prior, it was shown to achieve a convergence rate of order $O\left(N^{-\frac{1}{2}}\right)$. According to {\bf Example \ref{ex2}}, our Laguerre-based estimator with the choice of parameter $a=0$ generalized Laguerre function basis is appropriate. In such case, $\beta=7/3$ and thus, {our} minimax risk will be of order $O\left(N^{-\frac{3r^*}{3r^*+5}}\right)$.  Notice that our rate matches that of Tiwari and Zalkikar~(1990) for the special case of $r^*=5/3$, but it is faster than {theirs} for $r^* >5/3$ and slower for $r^*< 5/3$.
 \end{remark}
  \begin{remark}	  
{Ma and Balakrishnan~(2000) considered a kernel-based approach to the empirical Bayes estimation for the truncation parameter when the mixing distribution belongs to some general type of truncation parameter distributions, and were able to obtain a convergence rate of order $O\left(N^{-\frac{(\delta r-2)}{2r+1}}\right)$, where $\delta r >2$ and $0<\delta <2$, in their notation. This implies that if $r^*< \frac{\delta ((\beta\vee 1)+1)}{2(2-\delta)}$,our rate will be faster than theirs whenever $r < \frac{2[3r^*+(\beta\vee 1)+1]}{\delta(\beta\vee 1)-2r^*(2-\delta)}$ and slower otherwise. However, if $r^*> \frac{\delta ((\beta\vee 1)+1)}{2(2-\delta)}$, our rate will outperform theirs regardless of the value of $r$ in their notation}. 
  \end{remark}
  \begin{remark}	
 In the case of  mixing distributions with support $(0, \infty)$, such as the exponential, Rayleigh and the Weibull distributions, our estimator with different choices of the parameter $a$ of the generalized Laguerre function basis applies. Based on  {\bf Examples}  {\bf \ref{ex4}}, {\bf \ref{ex5}} and  {\bf \ref{ex6}}, $\beta=1$  and consequently, our posterior risk is of order $O\left(N^{-\frac{r^*}{r^*+1}}\right)$. Compared to Benhaddou and Pensky~(2013), in their wavelet-based procedure, and for $r^*=r$, their convergence rate for the general one-parameter exponential family of distributions is of order $O\left(N^{-\frac{2r}{2r+3}}\right)$, which is slightly slower than ours. 
  \end{remark}
  \begin{remark}	
 In the case of a $Beta (\alpha, \theta)$ mixing distribution with known $\alpha$ and unknown $\theta$ (or of the form $q(x\mid \theta)=\theta x^{\theta-1}$, for $\theta > 1$), which is an example of  mixing distributions that have positive compact supports, our estimator with choice of the parameter $a=2$ generalized Laguerre function basis is appropriate. Based on {\bf Example \ref{ex3}}, $\beta=1$  and therefore, our minimax risk is of order $O\left(N^{-\frac{r^*}{r^*+1}}\right)$. 
  \end{remark}
   \begin{remark}	
The optimal truncation level $M_o$ in \fr{R-min} depends on the smoothness parameter $r$ associated with the functions $p(x)$ and $\Psi(x)$ which is unknown, and therefore our estimator \fr{lagmat-hat} is not adaptive {with respect of $r$}. {Adaptivity may be achieved, for instance, using hard-thresholding on the estimated coefficients $(\widehat{{\pmb{{\Theta}}}^{\delta}}_M)_l$, $l=0, 1, 2, \cdots, M-1$, following the standard hard-thresholding procedure, or via Lepski's method (e.g., Lepski~(1991), Lepski, Mammen and Spokoiny~(1997)), for the optimal selection of $M$}. Providing completely data-driven procedure is beyond the scope of this work and we consider the parameter $r$ as if it was known. 
  \end{remark}

\section{Appendix} 

 {\bf Proof of Lemma \ref{lem:Bias}}. In order to prove the lemma, we follow three steps. The first step, we find a series expression for the inverse of matrix ${\bf A}_M$, ${\bf A}^{-1}_M$. In the second step, we evaluate the series expression for the matrix-vector product ${\bf A}^{-1}_M{\bf C}_M$ and then find an expression for $t_M(y)$ in terms of what we obtain from ${\bf A}^{-1}_M{\bf C}_M$, and in the last step we study the asymptotic behavior of $t_M(y)$ as $M \rightarrow \infty$. \\
   {\bf Finding Series expression for ${\bf A}^{-1}_M$}. Recall the series expansion of ${\bf A}_M$ in \fr{amser}, and let us find a series expansion for its inverse, ${\bf A}^{-1}_M$. In particular, suppose that ${\bf A}^{-1}_M$ has the form
 \beqn \label{aminvser}
{\bf A}^{-1}_M=\sum^r_{h=0}M^{-h}\frac{p^{(h)}(y)}{h!}\Xi^{(h)}_M+o(M^{-r}),
\eeqn   
 where  $\Xi^{(h)}_M$ are to be determined. Indeed, multiplying both sides of \fr{aminvser} and \fr{amser} yields
 \beqn \label{aminvprser}
{\bf A}_M{\bf A}^{-1}_M=\sum^{2r}_{i=0}M^{-i}\sum^{\min(i, r)}_{h=\max(0, i-r)}\frac{p^{(i-h)}(y)}{(i-h)!}{\bf \Phi}^{(i-h)}\Xi^{(h)}_M=I_M+o(M^{-r}).
\eeqn   
 Now, equating matrix coefficients for powers $i$ of $M^{-1}$ yields, for $i=0$, $i=1$ and $i=2$
 \beqns
 p(y){\bf \Phi}^{(0)}\Xi^{(0)}_M=I_M, \ \ \Xi^{(0)}_M=\frac{1}{p(y)}\left[{\bf \Phi}^{(0)}\right]^{-1}=\frac{1}{p(y)}I_M,
  \eeqns
 \beqns
p'(y){\bf \Phi}^{(1)}\Xi^{(0)}_M+p(y){\bf \Phi}^{(0)}\Xi^{(1)}_M=I_M, \ \ \Xi^{(1)}_M=\frac{-p'(y)}{p^2(y)}{\bf \Phi}^{(1)},
  \eeqns  
  and
\beqns
 \frac{p''(y)}{2}{\bf \Phi}^{(2)}\Xi^{(0)}_M+p'(y){\bf \Phi}^{(1)}\Xi^{(1)}_M+p(y){\bf \Phi}^{(0)}\Xi^{(2)}_M=I_M, \ \ \Xi^{(2)}_M=\frac{-p''(y)}{2p^2(y)}{\bf \Phi}^{(2)}+ \frac{(p'(y))^2}{p^3(y)}\left({\bf \Phi}^{(1)}\right)^2,
  \eeqns    
 respectively. The same way, for $i=3$, it can be shown that 
 \beqns
 \Xi^{(3)}_M= -\frac{p'''(y)}{6p^2(y)}{\bf \Phi}^{(3)}+\frac{p''(y)p'(y)}{2p^3(y)}{\bf \Phi}^{(2)}{\bf \Phi}^{(1)}+\frac{p''(y)p'(y)}{2p^3(y)}{\bf \Phi}^{(1)}{\bf \Phi}^{(2)}-\frac{(p'(y))^3}{p^4(y)}\left({\bf \Phi}^{(1)}\right)^3.
      \eeqns    
      Therefore, in general, for $i=j$, the matrices $\Xi^{(j)}_M$ in \fr{aminvser} will have the form 
      \be \label{xigen}
  \Xi^{(j)}_M=\sum \phi(j){\bf \Phi}^{(s_1)}{\bf \Phi}^{(s_2)}\cdots {\bf \Phi}^{(s_l)} ,         
  \ee
where ${\bf \Phi}^{(s)}$ are defined in \fr{amser} and $ \phi(j)$ depends on $p(y)$ and its derivatives up to the $j^{th}$ derivative and $\sum^{l}_{i=1}=j$.\\
 {\bf Finding series expression for ${\bf A^{-1}_M{\bf C}_M}$ and ${\bf t_M(y)}$}. To derive a series expression for ${\bf A^{-1}_M}{\bf C}_M$, multiply \fr{aminvser} and \fr{cmser} and use representation \fr{xigen} to obtain
\beqn
{\bf A}^{-1}_M{\bf C}_M&=&\sum^{2r}_{\nu=0}M^{-\nu}\sum^{\min(\nu, r)}_{h=\max(0, \nu-r)} \frac{\Psi^{(\nu-h)}(y)}{(\nu-h)!} \Xi^{(h)}_M\Lambda^{(\nu-h)}_M+o\left(M^{-r}\right)\nonumber\\
&=& \sum^{r}_{\nu=0}M^{-\nu}\sum^{\nu}_{h=0}\frac{\Psi^{(\nu-h)}(y)}{(\nu-h)!}\sum \phi(h){\bf \Phi}^{(s_1)}{\bf \Phi}^{(s_2)}\cdots {\bf \Phi}^{(s_l)}\Lambda^{(\nu-h)}_M+o\left(M^{-r}\right)
\eeqn
Consequently, 
\beqn
t_M(y)&=&\sum^{r}_{\nu=0}M^{-\nu}\sum^{\nu}_{h=0}\frac{\Psi^{(\nu-h)}(y)}{(\nu-h)!}\sum \phi(h)\sum^{M-1}_{l=0}\left({\bf \Phi}^{(s_1)}{\bf \Phi}^{(s_2)}\cdots {\bf \Phi}^{(s_l)}\Lambda^{(\nu-h)}_M\right)_l\varphi^{(a)}_l(y)+O\left(M^{-r}\right) \nonumber\\
&=& \frac{\Psi(y)}{p(y)}\sum^{M-1}_{l=0}\left(\Lambda^{(0)}_M\right)_l\varphi^{(a)}_l(y)+ \sum^{r}_{\nu=1}M^{-\nu}\sum^{\nu}_{h=0}\frac{\Psi^{(\nu-h)}(y)}{(\nu-h)!}\sum \phi(h)\sum^{M-1}_{l=0}\left({\bf \Phi}^{(s_1)}{\bf \Phi}^{(s_2)}\cdots {\bf \Phi}^{(s_l)}\Lambda^{(\nu-h)}_M\right)_l\varphi^{(a)}_l(y)\nonumber\\
&+& o\left(M^{-r}\right) \label{71}.
\eeqn
{\bf Asymptotic behavior for ${\bf t_M(y)}$}. Notice that, by {\bf Theorem $12$} of Muckenhoupt~(1970), which guarantees the mean convergence of the partial sum of the Laguerre polynomial series, one has 
\be \label{72}
\sum^{M-1}_{l=0}\left(\Lambda^{(0)}_M\right)_l\varphi^{(a)}_l(y)=\sum^{M-1}_{l=0}\langle1,\varphi^{(a)}_l(x) \rangle\varphi^{(a)}_l(y)  \rightarrow 1, as\ M \rightarrow \infty,
\ee
where $\langle f(x) , g(x) \rangle$ is the inner product between the functions $f$ and $g$. Also, all the terms in the second summation of \fr{71} will converge to zero, as $M \rightarrow \infty$.  We will only prove this for $\nu=1$. Denote the terms associated with $\nu=1$ and $\nu=j$ by $S_1$ and $S_j$, $j=2, 3, \cdots, r$, respectively. Indeed, for $\nu=1$, one has 
\beqn
S_1 &=&\frac{1}{M}\left[\frac{\Psi'(y)}{p(y)}\sum^{M-1}_{l=0}\langle (u-y), \varphi^{(a)}_l(x) \rangle\varphi^{(a)}_l(y)-\frac{p'(y)\Psi(y)}{p^2(y)}\sum^{M-1}_{l=0}\left({\bf \Phi}^{(1)}_M\Lambda^{(0)}_M\right)_l\varphi^{(a)}_l(y)\right]\nonumber\\
&= &\frac{1}{M}\left[\frac{\Psi'(y)}{p(y)}\sum^{M-1}_{l=0}\langle (u-y), \varphi^{(a)}_l(u) \rangle\varphi^{(a)}_l(y)\right]\nonumber\\
&-& \frac{1}{M} \left[\frac{p'(y)\Psi(y)}{p^2(y)}\sum^{M-1}_{l=0}\left(c_1\langle (u-y)\varphi^{(a)}_0(u), \varphi^{(a)}_l(u) \rangle+c_2\langle (u-y)\varphi^{(a)}_1(u), \varphi^{(a)}_l(u) \rangle \rangle\right)\varphi^{(a)}_l(y)\right]\nonumber\\
&+& \frac{1}{M} \left[\frac{p'(y)\Psi(y)}{p^2(y)}\sum^{M-1}_{l=0}\left(\cdots + c_{M}\langle (u-y)\varphi^{(a)}_{M}(u), \varphi^{(a)}_l(u) \rangle\right)\varphi^{(a)}_l(y)\right] \label{73}.
\eeqn
Now, denote $c_i=\int^{\infty}_0\varphi^{(a)}_{i-1}(u)du$, and notice that by \fr{73} and Theorem $12$ of Muckenhoupt~(1970), we have 
\be  \label{74} 
MS_1 \rightarrow \left[\frac{\Psi'(y)}{p(y)}(y-y)-\frac{p'(y)\Psi(y)}{p^2(y)}\left(c_1(y-y)\varphi^{(a)}_0(y) + c_2(y-y)\varphi^{(a)}_1(y) + \cdots + c_M(y-y)\varphi^{(a)}_{M-1}(y)\right) \right]=0.
\ee
The same way, we can show that 
\be \label{75}
M^{j}S_j \rightarrow 0,\ as\ M \rightarrow \infty.
\ee
Hence, by \fr{71}, \fr{72}, \fr{74} and \fr{75}, as $M \rightarrow \infty$, 
\beqn \label{t-asym}
t_M(y)=t(y)+o\left(M^{-(r\wedge s)}\right). 
\eeqn
To complete the proof, subtract $t(y)$ from both sides of \fr{t-asym} and square both sides. $\Box$\\
{\bf Proof of Lemma \ref{lem:var}}.   Notice that, with $A_{lk}$ defined in \fr{alk}, $\xi^{(l, k)}_i=\varphi_l(X_i)\varphi_k(X_i)-A_{lk}$, $i=1, 2, \cdots, N$, are zero mean independent and identically distributed having variance  
 \be
 \Var(\xi^{(l, k)}_i)\leq 2\int^{\infty}_0\left(\varphi^{(a)}_k(x) \varphi^{(a)}_l(x)\right)^2p(x)dx\leq 2\|p\|_{\infty}\min\{\|\varphi_l\|^2_{\infty}, \|\varphi_k\|^2_{\infty}\}.
 \ee    
 Therefore, by $(2.5)$ in Muckenhoupt~(1970), one has
  \be
\EE \| \widehat{\bf A}_M-{\bf A}_M\|_F^2=N^{-1} \sum_{k, l \leq M-1} Var(\xi^{(l, k)}_i)=O\left(N^{-1}M^2\right).
 \ee    
To complete the proof of \fr{var-a}, use the fact that the spectral norm of a matrix is less than or equal to its Frobenius norm. Similarly, 
\beqns
\EE\left|\widehat{A}_{lk}-A_{lk}\right|^4= O \left(N^{-4}\left[N\EE\left[\xi^{(l, k)}_i\right]^4+N(N-1)\EE^2\left[\xi^{(l, k)}_i\right]^2\right]\right)= O\left(N^{-2}\right), 
\eeqns 
and therefore,  since the matrix ${\bf A}_M$ is of size $M$, one has 
 \be
\EE \| \widehat{\bf A}_M-{\bf A}_M\|_{sp}^4=O\left(N^{-2}M^2\right).
 \ee    
 To prove \fr{var-c} and \fr{var-c2}, notice that, with $c_{k}$ defined in \fr{ck}, the quantities $\eta^{(k)}_i=U_k(X_i)-c_{k}$, $i=1, 2, \cdots, N$, are zero mean independent and identically distributed having variance  
 \be
 \Var(\eta^{(k)}_i)\leq 2\int^{\infty}_0\left(U_k(x)\right)^2p(x)dx\leq 2\|p\|_{\infty}\int^{\infty}_0U^2_k(x)dx.
 \ee    
Therefore, by  condition \fr{LRD}, one obtains 
 \be
\EE \| \widehat{{\bf C}}_M-{{\bf C}}_M\|^2=N^{-1} \sum_{k \leq M-1}  \Var(\eta^{(k)})\leq {\bf C}N^{-1}\sum_{k \leq M-1}(k^{\beta}\vee 1)=O\left(N^{-1}M^{\beta+1}\right).
 \ee    
 Similarly, 
\beqns
\EE\left|\widehat{c}_{k}-c_{k}\right|^4&=&O\left( N^{-4}\left[N\EE\left[\eta^{(k)}_i\right]^4+N(N-1)\EE^2\left[\eta^{(k)}_i\right]^2\right]\right)\\
&=& O\left(N^{-4}\left[N\int^{\infty}_{0}U^4_l(x)dx+N(N-1)\left[\int^{\infty}_0U^2_k(x)dx\right]^2\right]\right).
\eeqns 
To complete the proof of \fr{var-c2},  recall condition \fr{LRD} and note that the vectors above are of size $M$. $\Box$\\
The proof of {\bf Lemma \ref{lem:largd}} relies on the following version of Bernstein inequality. 
   \begin{lemma}(Bernstein Inequality). \label{lem:bernineq}
Let $Y_i$, $i=1, 2, \cdots, N$,  be independent and identically distributed with mean zero and finite variance $\sigma^2$, with $\|Y_i\| \leq \|Y\|_{\infty} < \infty$. Then,  
\be \label{prob-b}
\Pr\left(\left|N^{-1}\sum^N_{i=1}Y_i\right|>z\right) \leq 2 \exp\left\{-\frac{Nz^2}{2(\sigma^2+\|Y\|_{\infty}z/3)}\right\}.
\ee
\end{lemma}
{\bf Proof of Lemma \ref{lem:largd}}. 
Recall   \fr{ahat} and that the quantities $\xi^{(l, k)}_i=\varphi_l(X_i)\varphi_k(X_i)-A_{lk}$, $i=1, 2, \cdots, N$, are zero mean independent and identically distributed having variance $\sigma^2\leq 2\|p\|_{\infty}\|\varphi_l\|^2_{\infty}$, with $\|\xi^{(l, k)}\|_{\infty} < 2\|\varphi_l\|^2_{\infty}|\varphi_k\|^2_{\infty}$.  To apply Lemma \ref{lem:bernineq} take $z^2=\gamma^2\frac{\ln(N)}{N}$ to obtain  
\be \label{prob-b}
\Pr\left(\left|\widehat{A}_{lk}-A_{lk}\right| > z\right)= \Pr\left(\left|N^{-1}\sum^N_{i=1}\xi^{(l, k)}_i\right|>z\right) \leq 2 \exp\left\{-\frac{\gamma^2\ln(N)}{8\|\varphi_l\|^2_{\infty}\|p\|_{\infty}}\right\}.
\ee
To complete the proof, use result \fr{prob-b} and notice that  
\beqns
\Pr\left( \| \widehat{\bf A}_M-{\bf A}_M\|_2^2 > M^2 \gamma^2 N^{-1}\ln(N) \right)& \leq&2 \sum^{M-1}_{l=0}\sum^{M-1}_{k=0}\Pr\left(\left|\widehat{A}_{lk}-A_{lk}\right|^2 > z^2\right). \Box\label{devi-c}
\eeqns
 {\bf Proof of Lemma \ref{lem:R2}.} Let ${\bf A}_{\delta}={\bf A}+\delta {\bf I}$, where $\delta$ is a positive constant and ${\bf I}$ is the identity matrix. 
We will eventually choose the value of $\delta^2(N)$ where $\delta$ is an explicit function of $N$, as indicated in the previous sections, but in these proofs we will write using $\delta$ as an unknown constant.
Notice that
\beqn 
\|\widehat{{\pmb \Theta}}^{\delta}-{\pmb \Theta}\|\leq \|{\bf A}^{-1}\|_{sp}\|\widehat{{\bf C}}-{\bf C}\|+\|\widehat{\bf A}^{-1}_{\delta}-{\bf A}^{-1}\|_{sp}\|{\bf C}\|+\|\widehat{\bf A}^{-1}_{\delta}-{\bf A}^{-1}\|_{sp}\|\widehat{{\bf C}}-{\bf C}\|. \label{thethatdev1}
\eeqn
and
\beqn
\| \widehat{\bf A}^{-1}_{\delta}- {\bf A}^{-1}\|_{sp} \leq  \| \widehat{\bf A}^{-1}_{\delta}- {\bf A}^{-1}_{\delta}\|_{sp}+ \|{\bf A}^{-1}_{\delta}- {\bf A}^{-1}\|_{sp}. \label{thethatdev2}
\eeqn
Keep in mind that it is true that, for any nonsingular matrix ${\bf B}$, $\|{\bf B}^{-1}\|\geq \|{\bf B}\|^{-1}$. In addition, one can show that
 \be \label{ahata}
\widehat{\bf A}^{-1}_{\delta}-{\bf A}^{-1}_{\delta}=\widehat{\bf A}^{-1}_{\delta}(\widehat{\bf A}-{\bf A}){\bf A}^{-1}_{\delta}, \ \ {\bf A}^{-1}_{\delta}-{\bf A}^{-1}={\bf A}^{-1}_{\delta}({\bf A}-{\bf A}_{\delta}){\bf A}^{-1}.
\ee
  \begin{corollary}\label{Cor:bernineq}. For the nonsingular matrix ${\bf A}$ and positive scalar $\delta$, the following is true 
\be\label{adelta}
\|{\bf A}_{\delta}^{-1}\|_{sp}\leq \delta^{-1},\ \  \|\widehat{{\bf A}}_{\delta}^{-1}\|_{sp}\leq \delta^{-1}, \ \ \|{\bf A}_{\delta}^{-1}\|_{sp}\leq \|{\bf A}^{-1}\|_{sp}.
\ee
In addition, 
\be \label{vrs}
\|{\bf A}_{\delta}^{-1}-{\bf A}^{-1}\|_{sp}\leq \delta \|{\bf A}^{-1}\|_{sp}^2.
\ee
\end{corollary}
 {\bf Proof of Corollary \ref{Cor:bernineq}}. To  prove the three statements in \fr{adelta}, we start from the relation between the eigenvalue of a matrix and that of its inverse. Indeed,
\beqns
\|{\bf A}^{-1}_{\delta}\|_{sp}=\frac{1}{\|{\bf A}_{\delta}\|_{sp}}
=\frac{1}{\|{\bf A}\|_{sp}+\delta}\leq\delta^{-1}.
\eeqns 
A similar argument applies for $\|\widehat{{\bf A}}^{-1}_{\delta}\|_{sp}\leq\delta^{-1}$, and going in a different direction,
\beqns
\|{\bf A}^{-1}_{\delta}\|_{sp}=\frac{1}{\|{\bf A}\|_{sp}+\delta}\leq\frac{1}{\|{\bf A}\|_{sp}}= \|{\bf A}^{-1}\|_{sp}.\ \Box
\eeqns
Now, to prove \fr{vrs}, notice that by results \fr{ahata} and \fr{adelta} and the property of matrix norms, one has 
\beqns
\| {\bf A}_{\delta}^{-1}- {\bf A}^{-1}\|_{sp}=\| {\bf A}_{\delta}^{-1}( {\bf A}- {\bf A}_{\delta}) {\bf A}^{-1}\|_{sp}
\leq\|{\bf A}_{\delta}^{-1}\|_{sp}\|({\bf A}-{\bf A}_{\delta})\|_{sp}\|{\bf A}^{-1}\|_{sp}
\leq\delta\|{\bf A}^{-1}\|_{sp}^2. \ \Box
\eeqns
\\
Now, to find an upper bound for the first term of the right-hand side of \fr{thethatdev2}, introduce the sets $\Omega({\bf A})$ 
\be \label{omeg}
\Omega({\bf A})=\left\{\omega: \|\widehat{{\bf A}}-{\bf A}\| \geq 0.5 \|{\bf A}^{-1}\|^{-1}\right\}.
\ee
Notice that by inequality $||\widehat{{\bf A}}^{-1}_{\delta}-{\bf A}_{\delta}^{-1}|| \geq ||\widehat{{\bf A}}^{-1}_{\delta}||-||{\bf A}_{\delta}^{-1}||$,  it is easy to show that
\beqn \label{ainvdel}
\|\widehat{{\bf A}}^{-1}_{\delta}\|_{sp} \leq 2 \|{\bf A}^{-1}\|_{sp}.
\eeqn
Therefore, by \fr{ainvdel} and inequality $||\widehat{{\bf A}}^{-1}_{\delta} -{\bf A}^{-1}_{\delta}||_{sp} \leq ||{\bf A}^{-1}_{\delta}||_{sp}$ combined with \fr{adelta}, the first term in the right-hand side of \fr{thethatdev2} is such that 
\beqn
\|\widehat{{\bf A}}^{-1}_{\delta} -{\bf A}^{-1}_{\delta}\|_{sp} &\leq& 2\|{\bf A}^{-1}\|_{sp}^2\|\widehat{{\bf A}}-{\bf A}\|_{sp}\II(\Omega^c({\bf A})) + \|\widehat{{\bf A}}^{-1}_{\delta} -{\bf A}^{-1}_{\delta}\|_{sp}\II(\Omega({\bf A}))\nonumber\\
&\leq& 2\|{\bf A}^{-1}\|_{sp}^2\|\widehat{{\bf A}}-{\bf A}\|_{sp} + 2\delta^{-1}\II(\Omega({\bf A})).\label{varalep}
\eeqn
Putting all this together, we reconstruct (\ref{thethatdev1}) into
\be
\|\widehat{{\pmb \Theta}}^{\delta}-{\pmb \Theta}\|\leq \|{\bf A}^{-1}\|_{sp}\|\widehat{{\bf C}}-{\bf C}\|+\left(\| \widehat{{\bf A}}^{-1}_{\delta}- {\bf A}^{-1}_{\delta}\|_{sp}+ \|{{\bf A}}^{-1}_{\delta}- {\bf A}^{-1}\|_{sp}\right)\left(\|{\bf C}\|+\|\widehat{{\bf C}}-{\bf C}\|\right). \label{thethatdev3}\\
\ee
And then we use (\ref{vrs}) and (\ref{varalep}) to show
\beqn
\|\widehat{{\pmb \Theta}}^{\delta}-{\pmb \Theta}\| &\leq& \|{\bf A}^{-1}\|_{sp}\|\widehat{{\bf C}}-{\bf C}\| \nonumber\\
&+& \left(2\|{\bf A}^{-1}\|_{sp}^2\|\widehat{{\bf A}}-{\bf A}\|_{sp} + 2\delta^{-1}\II(\Omega({\bf A}))\right)\left(\|{\bf C}\|+\|\widehat{{\bf C}}-{\bf C}\|\right)\nonumber\\
&+&  \left(\delta \|{\bf A}^{-1}\|_{sp}^2\right)\left(\|{\bf C}\|+\|\widehat{{\bf C}}-{\bf C}\|\right). \label{thethatdev3}
\eeqn
When squaring, we will use the fact that, for any $\bf X$ and $\bf Y$,
\beqns
\left(\|\bf X\|+\|\bf Y\|\right)^2 \leq 2\left(\|\bf X\|\right)^2 + 2\left(\|\bf Y\|\right)^2,
\eeqns
to find that
\beqn
\|\widehat{{\pmb \Theta}}^{\delta}-{\pmb \Theta}\|^2 &\leq& 2\|{\bf A}^{-1}\|^2_{sp}\|\widehat{{\bf C}}-{\bf C}\|^2 \nonumber\\
&+& 4\|{\bf A}^{-1}\|^4_{sp}\|\widehat{{\bf A}}-{\bf A}\|^2_{sp}\|{\bf C}\|^2+ 4\|{\bf A}^{-1}\|^4_{sp}\|\widehat{{\bf A}}-{\bf A}\|^2_{sp}\|\widehat{{\bf C}}-{\bf C}\|^2\nonumber\\
&+& 4\delta^{-2}\II(\Omega({\bf A}))\|{\bf C}\|^2+ 4\delta^{-2}\II(\Omega({\bf A}))\|\widehat{{\bf C}}-{\bf C}\|^2\nonumber\\
&+&  2\delta^2 \|{\bf A}^{-1}\|^4_{sp}\|{\bf C}\|^2+2\delta^2\|{\bf A}^{-1}\|^4_{sp}\|\widehat{{\bf C}}-{\bf C}\|^2.          \label{thethatdev4}
\eeqn
Then, we take the expectation and apply the Cauchy-Schwarz Inequality for expectation, that for any random quantities $\bf X$ and $\bf Y$
\beqns
|\EE {\bf {XY}}| \leq \sqrt{\EE \bf X^2 \EE \bf Y^2},
\eeqns
 We thus find
\beqn
\EE\|\widehat{{\pmb \Theta}}^{\delta}-{\pmb \Theta}\|^2&=&O\left( \EE\|\widehat{{\bf C}}-{\bf C}\|^2+\EE\|\widehat{{\bf A}}-{\bf A}\|_{sp}^2+\sqrt{\EE\|\widehat{{\bf A}}-{\bf A}\|_{sp}^4\EE\|\widehat{{\bf C}}-{\bf C}\|^4}+\frac{1}{\delta^2}\Pr(\Omega({\bf A}))\right)\nonumber\\
&+& O\left(\frac{1}{\delta^2}\sqrt{\EE\|\widehat{{\bf C}}-{\bf C}\|^4\Pr(\Omega({\bf A}))}+\delta^2+ \delta^2\EE\|\widehat{{\bf C}}-{\bf C}\|^2\right) \label{thethatdev5}
\eeqn
We use the earlier-proven expectation values
\beqn
\EE \| \widehat{\bf A}_M-{\bf A}_M\|_{sp}^2 &\leq& O\left(N^{-1}M^2\right).\\
\EE \| \widehat{\bf A}_M-{\bf A}_M\|_{sp}^4&=&O\left(N^{-2}M^2\right).\\
\EE \| \widehat{{\bf C}}_M-{{\bf C}}_M\|^2&=&O\left(N^{-1}M^{\beta+1}\right).\\
\EE\|\widehat{{\bf C}}_M-{{\bf C}}_M\|^4 &=&O\left(1\vee M^{2\beta+1}N^{-2}\right),
\eeqn
the last of which is found through
\beqn
\EE\|\widehat{{\bf C}}_M-{{\bf C}}_M\|^4 = O\left(M N^{-4}\max_{l\leq M-1}\left[N\int^{\infty}_{0}U^4_l(x)dx+N(N-1)\left[\int^{\infty}_0U^2_k(x)dx\right]^2\right]\right),
\eeqn
combined with the fact that $M\max_{l\leq M-1}\int^{\infty}_0U^4_l(x)dx=o(N^3)$ and  $\int^{\infty}_0U^2_{M-1}(x)dx=O(M^{\beta})$. From this set of expectations, and now defining $\delta^2(N)=M^2N^{-1}$, we find that \fr{thethatdev5} yields
\beqn
\EE\|\widehat{{\pmb \Theta}}^{\delta}-{\pmb \Theta}\|^2 &=& O\left( \frac{M^{\beta+1}}{N} +\frac{M^2}{N} + \frac{N}{M^2}\Pr(\Omega({\bf A}))+\frac{N}{M^2}\frac{M^{\beta+1/2}}{N}\sqrt{\Pr(\Omega({\bf A}))}\right) \label{thethatdev7}
\eeqn
Finally, taking $\gamma^2=\frac{N\|{\bf A}^{-1}\|^2}{4M^2\ln(N)}$ in \fr{thethatdev7} completes the proof.  $\Box$\\
 {\bf Proof of Theorem \ref{th:upperbds-2}}.   Combining \fr{var-bias} and \fr{R-min} in \fr{rn12} and plugging in \fr{R-min} completes the proof. $\Box$

\end{document}